\documentclass[a4, 12pt]{article}
\usepackage{amsmath}

\makeatletter
\renewcommand*\env@matrix[1][*\c@MaxMatrixCols c]{%
  \hskip -\arraycolsep
  \let\@ifnextchar\new@ifnextchar
  \array{#1}}
\makeatother
\usepackage{mathptmx}
\usepackage{amssymb}
\usepackage{pdfpages}
\usepackage{graphics}
\usepackage{tikz}
\usepackage[margin=1in]{geometry}
\usepackage{setspace}

\newtheorem{theorem}{Theorem}

\numberwithin{subcase}{case}
\newtheorem{lemma}{Lemma}

\newcommand{\qed}{\hfill$\Box$}

\begin{document}
\title{On Kirchhoff index and number of spanning trees of linear pentagonal cylinder and M\"{o}bius chain graph}
\author{Md. Abdus Sahir\thanks{Email: abdussahir@gmail.com} ~and Sk. Md. Abu Nayeem\thanks{Corresponding author. Email: nayeem.math@aliah.ac.in}\\{\small Department of Mathematics and Statistics, Aliah University, Kolkata -- 700 160, India.}}
\date{}
\maketitle

\noindent
\begin{abstract} 
In this paper, we derive closed-form formulas for Kirchhoff index and Wiener index of linear pentagonal cylinder graph and linear pentagonal M\"{o}bius chain graph. We also obtain explicit formulas for finding total number of spanning trees for both the graphs.

\medskip



\end{abstract}

\onehalfspacing
\section{Introduction}
Let $G=(V,E)$ be a connected molecular graph having vertex set $V$ and edge set $E$. A topological index of $G$ is a numerical quantity involving different graph parameters such as number of vertices, number of edges, degree, eccentricity, distance between vertices, etc. Winner index, one of the oldest topological indices, was introduced by Harold Wiener \cite{Wie47} and was defined as $W(G)=\sum\limits_{i<j} d_{ij}$, where $d_{ij}$ is the length of the shortest path between the vertices $i$ and $j$. Motivated by the idea of Wiener index, the idea of resistance distance and Kirchhoff index was introduced by Klein and Randi\'{c} \cite{Kle93}. Kirchhoff index, initially known as resistance index, was defined as $Kf(G)=\sum\limits_{i<j} r_{ij}$, where $r_{ij}$ is the effective resistance between vertex $i$ and vertex $j$ calculated using Ohm's law considering all the edges of $G$ as unit resistors. Klein and Randi\'{c} also proved that for any vertex pair $i,j$ in a graph $G$, $r_{ij} \leq d_{ij}$ and $Kf(G)\leq W(G)$ with equality holds if and only if $G$ is a tree. 

In the last few decades, researchers have focused on various topological indices such as Wiener index, Randi\'{c} index \cite{Ran75}, Kirchhoff index, Gutman index \cite{Gut94}, Estrada index \cite{Est00}, Zagreb index \cite{Gut72} etc. Especially Kirchhoff index is seeking a lot of attention of researchers as it has wide applications in physics, chemistry, graph theory and various related subjects. Readers are referred to \cite{Est10, Mad13, Wan10, Wan11, Xia03, You13, Zho08} for some recent works. Many researchers have concentrated on finding the Kirchhoff index and the number of spanning trees for many interesting graphs such as linear hexagonal chain \cite{Yan08}, linear pentagonal chain \cite{Wan10}, M\"{o}bius hexagonal chain \cite{Wan11}, periodic linear chain \cite{Car15}, crossed hexagonal chain \cite{Pan18}, linear octagonal chain \cite{Zhu20}, M\"{o}bius/ cylinder octagonal chain \cite{Liu22} and many others \cite{Gen21, Liu19, Lui21, Pen17}. Normalized Laplacian spectrum and the number of spanning trees of linear pentagonal chains have been obtained by He et al. \cite{He18}. Although the Kirchhoff index of linear pentagonal chain was found long back in 2010 \cite{Wan10}, but to the best of our knowledge, Kirchhoff indices for linear pentagonal cylinder and M\"{o}bius chains have not been obtained so far. In the present paper, we aim to obtain those.
 
The Laplacian matrix $L(G)=(l_{ij})$ of a simple connected graph $G$ is defined as,
\[   
l_{ij}= 
     \begin{cases}
       d_i, &\mbox{if } i=j\\
      -1,  & \mbox{if } i\sim j\\
       0, &\mbox{otherwise,}\\
     \end{cases}
\]
where $d_i$ is the degree of the vertex $i$. 

Since $L(G)$ is a symmetric matrix, all of its eigenvalues are real. Moreover all are non-negative, i.e., if the eigenvalues $\lambda_i$,s $(i=1,2,\ldots,n)$ are indexed in the increasing order of their values, $0=\lambda_1\leq \lambda_2 \leq \cdots \leq \lambda_n$. Since $G$ is connected, $\lambda_1=0$ is a simple eigenvalue. 

Gutman and Mohar \cite{Gut96} and Zhu et al. \cite{Zhu96} obtained the following lemma. 
\begin{lemma}  \cite{Gut96,Zhu96} For a connected graph $G$ with $n$-vertices, $n\geq 2$,
$$Kf(G)=n\sum\limits_{k=2}^{n}\frac{1}{\lambda_k}\cdot$$
\end{lemma}

\begin{figure}
\begin{center}
\begin{tikzpicture}[place/.style={thick, circle,draw=black!50,fill=black!20,inner sep=0pt, minimum size = 6 mm}, transform shape]
\node[place] (1) {\footnotesize 1} ;
\node (2)[place] [below of = 1] {\footnotesize $1'$} edge[thick] (1);
\node (3)[place] [xshift=1cm, above of = 1] {\footnotesize 2} edge[thick] (1);
\node (4)[place] [xshift=1cm, below of = 2] {\footnotesize $2'$} edge[thick] (2);
\node (5)[place] [below of = 3, node distance=1.5cm] {\footnotesize $\tilde{1}$} edge[thick] (3) edge[thick] (4);
\node (6)[place] [xshift=1cm, below of = 3] {\footnotesize 3} edge[thick] (3);
\node (7)[place] [below of = 6] {\footnotesize $3'$} edge[thick] (6) edge[thick] (4);
\node (8)[place] [xshift=1cm, above of = 6] {\footnotesize 4} edge[thick] (6);
\node (9)[place] [below of = 8, node distance=1.5cm] {\footnotesize $\tilde{2}$} edge[thick] (8) ;
\node (10)[place] [xshift=1cm, below of = 7] {\footnotesize $4'$} edge[thick] (9) edge[thick] (7);
\node (11)[place] [xshift=1cm, below of = 8] {\footnotesize 5} edge[thick] (8);
\node (12)[place] [below of = 11] {\footnotesize $5'$} edge[thick] (10) edge[thick] (11);
\node (13)[place] [right of = 11, node distance=1.5cm] {}  edge [thick, dotted] (11);
\node (14)[place] [below of = 13] {} edge[thick] (13) edge [thick, dotted] (12);
\node (15)[place] [xshift=1cm, above of = 13] {} edge[thick] (13);
\node (16)[place] [xshift=1cm, below of = 14] {} edge[thick] (14);
\node (17)[place] [below of = 15, node distance=1.5cm] {} edge[thick] (15) edge[thick] (16);
\node (18)[place] [xshift=1cm, below of = 15] {\footnotesize \tiny $2n-1$} edge[thick] (15);
\node (19)[place] [below of = 18] {\footnotesize \tiny $2n-1'$} edge[thick] (18) edge[thick] (16);
\node (20)[place] [xshift=1cm, above of = 18] {\footnotesize $2n$} edge[thick] (18)  edge[thick, out=125, in=125,out looseness=0.4, in looseness=.8] (1);
\node (21)[place] [below of = 20, node distance=1.5cm] {\footnotesize $\tilde{n}$} edge[thick] (20);
\node (22)[place] [xshift=1cm, below of = 19] {\footnotesize $2n'$} edge[thick] (21) edge[thick] (19) edge[thick, out=-150, in=-125, out looseness=0.4, in looseness=.8] (2);
\end{tikzpicture}
\caption{\label{fig1} Linear pentagonal cylinder graph $P_n.$}
\end{center}
\end{figure}
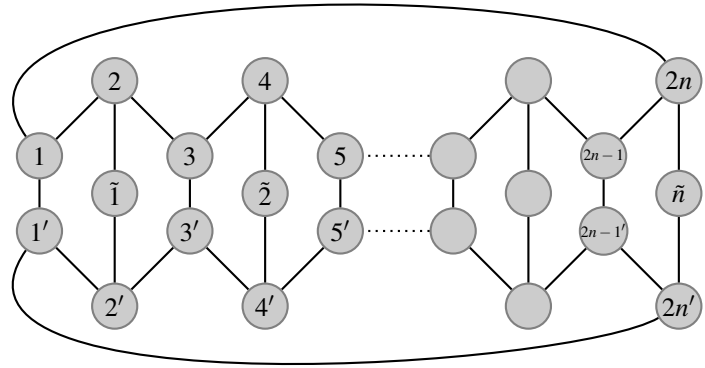

\begin{figure}
\begin{center}
\begin{tikzpicture}[place/.style={thick, circle,draw=black!50,fill=black!20,inner sep=0pt, minimum size = 6 mm}, transform shape]
\node[place] (1) {\footnotesize 1} ;
\node (2)[place] [below of = 1] {\footnotesize $1'$} edge[thick] (1);
\node (3)[place] [xshift=1cm, above of = 1] {\footnotesize 2} edge[thick] (1);
\node (4)[place] [xshift=1cm, below of = 2] {\footnotesize $2'$} edge[thick] (2);
\node (5)[place] [below of = 3, node distance=1.5cm] {\footnotesize $\tilde{1}$} edge[thick] (3) edge[thick] (4);
\node (6)[place] [xshift=1cm, below of = 3] {\footnotesize 3} edge[thick] (3);
\node (7)[place] [below of = 6] {\footnotesize $3'$} edge[thick] (6) edge[thick] (4);
\node (8)[place] [xshift=1cm, above of = 6] {\footnotesize 4} edge[thick] (6);
\node (9)[place] [below of = 8, node distance=1.5cm] {\footnotesize $\tilde{2}$} edge[thick] (8) ;
\node (10)[place] [xshift=1cm, below of = 7] {\footnotesize $4'$} edge[thick] (9) edge[thick] (7);
\node (11)[place] [xshift=1cm, below of = 8] {\footnotesize 5} edge[thick] (8);
\node (12)[place] [below of = 11] {\footnotesize $5'$} edge[thick] (10) edge[thick] (11);
\node (13)[place] [right of = 11, node distance=1.5cm] {}  edge [thick, dotted] (11);
\node (14)[place] [below of = 13] {} edge[thick] (13) edge [thick, dotted] (12);
\node (15)[place] [xshift=1cm, above of = 13] {} edge[thick] (13);
\node (16)[place] [xshift=1cm, below of = 14] {} edge[thick] (14);
\node (17)[place] [below of = 15, node distance=1.5cm] {} edge[thick] (15) edge[thick] (16);
\node (18)[place] [xshift=1cm, below of = 15] {\footnotesize \tiny $2n-1$} edge[thick] (15);
\node (19)[place] [below of = 18] {\footnotesize \tiny $2n-1'$} edge[thick] (18) edge[thick] (16);
\node (20)[place] [xshift=1cm, above of = 18] {\footnotesize $2n$} edge[thick] (18)  edge[thick, out=125, in=125,out looseness=0.4, in looseness=1.3] (2);
\node (21)[place] [below of = 20, node distance=1.5cm] {\footnotesize $\tilde{n}$} edge[thick] (20);
\node (22)[place] [xshift=1cm, below of = 19] {\footnotesize $2n'$} edge[thick] (21) edge[thick] (19) edge[thick, out=-150, in=-125, out looseness=0.4, in looseness=1.3] (1);
\end{tikzpicture}
\caption{\label{fig2} Pentagonal M\"{o}bius chain graph $P'_n.$}
\end{center}
\end{figure}
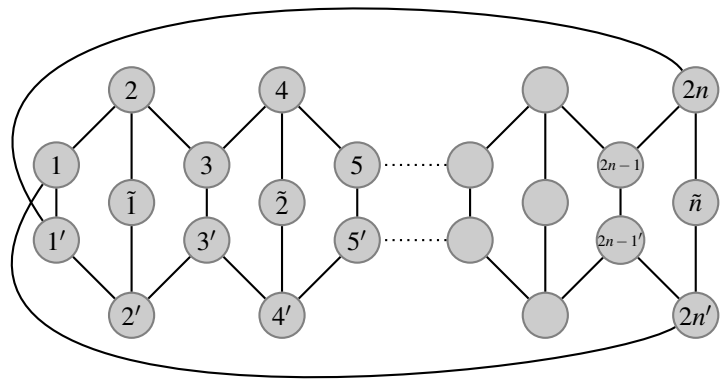

Like Wiener index, Kirchhoff index also gives description of the underlying structure of a molecular graph \cite{Xia03}. Obtaining closed-form formulae for Kirchhoff index of general graphs is not straight forward, but we can derive closed-form formulae for some special classes of graphs like cycles \cite{Kle95}, complete graphs \cite{Luk99}, circulant graph, etc. 

In this paper, we derive closed-form formulas for Kirchhoff index and Wiener index of linear pentagonal cylinder graph $P_n$ (Figure \ref{fig1}) and pentagonal M\"{o}bius chain graph $P'_n$ (Figure \ref{fig2}) on $5n$ ($n\geq 2$) vertices. Also we present the formulas for the total number of spanning trees for those graphs. 

\section{Preliminaries} Let $G$ be a graph with vertex set $V(G)$. A permutation $\pi$ of $V(G)$ is called an automorphism if $u$ and $v$ are adjacent in $G$ if and only if $\pi(u)$ and $\pi(v)$ are also adjacent in $G$. Suppose $V_0=\{\tilde{1},\tilde{2},\ldots,\tilde{p}\}$, $V_1=\{1,2,\ldots,q\}$ and $V_2=\{1',2',\ldots,q'\}$ are the vertex partition for the automorphism $\pi$ such that $\pi(\tilde{i})=\tilde{i}\mbox{ for all } \tilde{i}\in V_0, \pi(i)=i'\mbox{ for all }i\in V_1$ and $\pi(i')=i\mbox{ for all }i'\in V_2$. It is easy to follow that $\pi$ can be decomposed as product of disjoint $1$-cycles and transpositions, i.e.,
$$\pi=(\tilde{1})(\tilde{2})\cdots(\tilde{p})(1,1')(2,2')\cdots(q,q'),$$ where $p+2q=|V(G)|$. Then by suitable arrangement of vertices, the Laplacian matrix $L(G)$ of $G$ can be expressed into the block matrix form --
$$L(G)=\begin{bmatrix}
L_{V_0V_0} & L_{V_0V_1} & L_{V_0V_2}\\
L_{V_1V_0} & L_{V_1V_1} & L_{V_1V_2}\\
L_{V_2V_0} & L_{V_2V_1} & L_{V_2V_2}
\end{bmatrix}$$ where the submatrix $L_{V_rV_s}$ correspond to the vertices of $V_r$ and $V_s,r,s=0,1,2$ respectively. 

Let $$L_A(G)=\begin{bmatrix}L_{V_0V_0} & \sqrt{2}L_{V_0V_1} \\
\sqrt{2}L_{V_1V_0} & L_{V_1V_1} + L_{V_1V_2}\\
\end{bmatrix}$$ and $$L_S(G)=L_{V_1V_1} - L_{V_1V_2}.$$
Yang and Yu \cite{Yan85} and many others like Yang and Zhang \cite{Yan08} have used the Laplacian decomposition formula to find the Kirchhoff indices of certain classes of graphs where some automorphisms are found. We describe it in the form of the following lemma.  
\begin{lemma}\label{Yan}
The characteristic polynomial of $L(G)$ is equal to the product of that of $L_A(G)$ and $L_S(G)$, i.e., $$\det(L(G)-\lambda I)=\det(L_A(G)-\lambda I)\cdot\det(L_S(G)-\lambda I).$$
\end{lemma}

Let $\lambda_i$ $(i=1,2,\ldots, 3n)$ and $\mu_j$ $(j=1,2,\ldots,2n)$ are the eigenvalues of $L_A(G)$ and $L_S(G)$ arranged in ascending order of their values. Then by Lemma \ref{Yan}, the spectrum of $L(G)$ is given by $$\{0=\lambda_1\le\lambda_2\le\cdots\le\lambda_{3n}\}\bigcup \{\mu_1\le\mu_2\le\cdots\le\mu_{2n}\}.$$

To avoid confusion, we denote $L_A(P_n)$ and $L_S(P_n)$ by $L_A$ and $L_S$ respectively and $L_A(P_n')$ and $L_S(P_n')$ by $L'_A$ and $L'_S$ respectively. Also we denote the block matrices constituting the Laplacian matrix of linear pentagonal cylinder graph $P_n$ by $L_{V_rV_s},r,s=0,1,2$ and those of pentagonal M\"{o}bius chain graph $P'_n$ by $L'_{V_rV_s},r,s=0,1,2$ respectively. Then,

$L_{V_0V_0}=L'_{V_0V_0}=\begin{bmatrix}
2 & 0& \cdots &0\\
0 & 2 & \cdots & 0\\
\vdots & \vdots & \ddots & \vdots\\
0 & 0& \cdots &2
\end{bmatrix}_{n\times n}$,
$L_{V_0V_1}=L'_{V_0V_1}=\begin{bmatrix}
0 & -1&0&0& \cdots &0\\
0 & 0&0&-1& \cdots &0\\
\vdots & \vdots & \vdots&\vdots  &\ddots &  \vdots\\
0 & 0&0&0& \cdots &-1\\
\end{bmatrix}_{n\times 2n}$

$L_{V_1V_1}=\begin{bmatrix}
3 &-1&0& \cdots&0&-1\\
-1 &3&-1& \cdots&0&0\\
0 & -1&3& \cdots &0&0\\
\vdots & \vdots & \vdots&\ddots   & \vdots& \vdots\\
0 &0&0& \cdots&3&-1\\
-1 &0&0& \cdots&-1&3
\end{bmatrix}_{2n\times 2n}$,
$L_{V_1V_2}=\begin{bmatrix}
-1 &0&0& \cdots&0&0\\
0 &0&0& \cdots&0&0\\
0 & 0&-1& \cdots &0&0\\
\vdots & \vdots & \vdots&\ddots   & \vdots& \vdots\\
0 &0&0& \cdots&-1&0\\
0 &0&0& \cdots&0&0
\end{bmatrix}_{2n\times 2n}$,
and $L'_{V_1V_1}=\begin{bmatrix}
3 &-1&0& \cdots&0&0\\
-1 &3&-1& \cdots&0&0\\
0 & -1&3& \cdots &0&0\\
\vdots & \vdots & \vdots&\ddots   & \vdots& \vdots\\
0 &0&0& \cdots&3&-1\\
0 &0&0& \cdots&-1&3
\end{bmatrix}_{2n\times 2n}$,
$L'_{V_1V_2}=\begin{bmatrix}
-1 &0&0& \cdots&0&-1\\
0 &0&0& \cdots&0&0\\
0 & 0&-1& \cdots &0&0\\
\vdots & \vdots & \vdots&\ddots   & \vdots& \vdots\\
0 &0&0& \cdots&-1&0\\
-1 &0&0& \cdots&0&0
\end{bmatrix}_{2n\times 2n}$.\\
So, 
$$L_A=L_A'=\begin{bmatrix}[cccc|cccccc]
2 & 0 & \cdots & 0  & 0 & {-\sqrt{2}} & 0 & 0 & \cdots & 0\\
0 &2&\cdots&0  &0 &0&0&-\sqrt{2}&\cdots&0 \\
\vdots&\vdots&\ddots&\vdots&\vdots&\vdots&\vdots&\vdots&\ddots&\vdots\\
0 & 0 & \cdots & 2  & 0 & 0 & 0 & 0 & \cdots & {-\sqrt{2}}\\\hline

0 & 0 & \cdots & 0           & 2& -1 & 0 & 0 & \cdots & -1\\
{-\sqrt{2}} & 0 & \cdots & 0    & -1& 3 & -1 & 0 & \cdots & 0\\
0 & 0 & \cdots & 0              & 0& -1 & 2 & -1 & \cdots & 0\\
0 & {-\sqrt{2}}& \cdots & 0    & 0& 0& -1 & 3 & \cdots & 0\\
\vdots&\vdots&\ddots&\vdots&\vdots&\vdots&\vdots&\vdots&\ddots&\vdots\\
0 & 0& \cdots & {-\sqrt{2}}   & -1& 0& 0 & 0 & \cdots & 3\\
\end{bmatrix}_{3n\times 3n},$$

$$L_S=\begin{bmatrix}
4 &-1&0&0& \cdots&0&0&-1\\
-1&3 &-1&0& \cdots&0&0&0\\
0 &-1&4&-1& \cdots&0&0&0\\
0&0&-1&3 & \cdots&0&0&0\\
\vdots & \vdots & \vdots   & \vdots& \ddots& \vdots & \vdots & \vdots \\
0&0&0&0 & \cdots&3&-1&0\\
0&0&0&0 & \cdots&-1&4&-1\\
-1&0&0&0 & \cdots&0&-1&3\\
\end{bmatrix}_{2n\times 2n},$$

\begin{eqnarray*}\mbox{and }L_S'=\begin{bmatrix}
4 &-1&0&0& \cdots&0&0&1\\
-1&3 &-1&0& \cdots&0&0&0\\
0 &-1&4&-1& \cdots&0&0&0\\
0&0&-1&3 & \cdots&0&0&0\\
\vdots & \vdots &\vdots   &\vdots & \ddots&  \vdots & \vdots & \vdots \\
0&0&0&0 & \cdots&3&-1&0\\
0&0&0&0 & \cdots&-1&4&-1\\
1&0&0&0 & \cdots&0&-1&3
\end{bmatrix}_{2n\times 2n}.\end{eqnarray*}

For our convenience, we denote the matrix $$\begin{bmatrix}
4 &-1&0&0& \cdots&0&0&0\\
-1&3 &-1&0& \cdots&0&0&0\\
0 &-1&4&-1& \cdots&0&0&0\\
0&0&-1&3 & \cdots&0&0&0\\
\vdots & \vdots &\vdots   &\vdots & \ddots&  \vdots & \vdots & \vdots \\
0&0&0&0 & \cdots&3&-1&0\\
0&0&0&0 & \cdots&-1&4&-1\\
0&0&0&0 & \cdots&0&-1&3
\end{bmatrix}_{2n\times 2n}$$ by $L^0_S$, so that $L_S=L^0_S-e_1e_{2n}^T$ and $L_S'=L^0_S+e_1e_{2n}^T$ where $e_i$ is the unit column vector of compatible size with all its components 0 except the $i$th component which has the value 1.

In this paper, we shall use the following lemma, known as the matrix-determinant lemma to compute the determinant of a matrix with rank one perturbation if the determinant of the original matrix is known. 
\begin{lemma}\label{Matdet}
Let $M$ be an $n\times n$ matrix. Then $\det(M+uv^T)=\det(M)+v^T\operatorname{adj}(M)u,$ where $u,v$ are $n\times 1$ column vectors.
\end{lemma}

\section{Kirchhoff index of pentagonal cylinder and M\"{o}bius chain}
From Lemma \ref{Yan}, we have that the Kirchhoff index of linear pentagonal cylinder $P_n$ is $$Kf(P_n)=5n\left(\sum\limits_{i=2}^{3n} \frac{1}{\rho_i}+\sum\limits_{j=1}^{2n} \frac{1}{\mu_j}\right), n\geq 2$$ where $\rho_i,i=1,2,\ldots,3n$ and $\mu_j,j=1,2,\ldots,2n$ are the eigenvalues of $L_A$ and $L_S$ respectively. 

Let \begin{eqnarray}\det(xI_{3n}-L_A)=x^{3n}+\alpha_1x^{3n-1}+\cdots+\alpha_{3n-2}x^{2}+\alpha_{3n-1}x, (\mbox{since } \rho_1=0)\label{char1}\end{eqnarray} and \begin{eqnarray}\det(xI_{2n}-L_S)=x^{2n}+\beta_1x^{2n-1}+\cdots+\beta_{2n-2}x^{2}+\beta_{2n-1}x+\beta_{2n}.\label{char2}\end{eqnarray}

From Vieta's formula, we have $\sum\limits_{i=2}^{3n} \frac{1}{\rho_i}=-\frac{\alpha_{3n-2}}{\alpha_{3n-1}}$ and $\sum\limits_{j=1}^{2n} \frac{1}{\mu_i}=\frac{\beta_{2n-1}}{\beta_{2n}}=\frac{\beta_{2n-1}}{\det(L_S)}\cdot$ 

Hence $$Kf(P_n)=5n\left(-\frac{\alpha_{3n-2}}{\alpha_{3n-1}}+\frac{\beta_{2n-1}}{\det(L_S)}\right)\cdot$$

By similar argument, \begin{eqnarray*}Kf(P_n')&=&5n\left(\sum\limits_{i=2}^{3n} \frac{1}{\rho_i}+\sum\limits_{j=1}^{2n} \frac{1}{\mu_j'}\right)\\
&=&5n\left(-\frac{\alpha_{3n-2}}{\alpha_{3n-1}}+\frac{\beta'_{2n-1}}{\det(L_S')}\right),
\end{eqnarray*}
where $\mu_j',j=1,2,\ldots,2n$ are the eigenvalues of $L_S'$ and $\beta'_{2n-1}$ is the coefficient of the first degree term in the characteristic polyomial of $L_S'$ respectively.

\begin{lemma}\label{det}
Let $R_n=\begin{bmatrix}
-2 &1&0& \cdots&0&0\\
1 &-2&1& \cdots&0&0\\
0 & 1&-2& \cdots &0&0\\
\vdots & \vdots & \vdots&\ddots   & \vdots& \vdots\\
0 &0&0& \cdots&-2&1\\
0 &0&0& \cdots&1&-2
\end{bmatrix}_{n\times n},$ then $\det(R_n)=(-1)^n(1+n).$
\end{lemma}
\noindent\textit{Proof.} Here $\det(R_1)=-2,~\det(R_2)=3$ and $\det(R_n)=-2\det(R_{n-1})-\det(R_{n-2}),~n\geq3.$ Solving the recurrence relation, we get $\det(R_n)=(-1)^n(1+n).$ \qed 

\begin{lemma}\label{det1.1}
Let $$R_{n,m}= 
\begin{bmatrix}[ccccc|c|cccc]-2 &1&0&\cdots&0&0& 0&\cdots&0&0\\
1&-2 &1&\cdots&0&0& 0&\cdots&0&0\\
0 &1&-2&\cdots&0&0& 0&\cdots&0&0\\
\vdots & \vdots & \vdots& \ddots&\vdots   & \vdots&  \vdots& \ddots& \vdots & \vdots  \\
0&0&0&\cdots&-2 &1& 0&\cdots&0&0\\\hline 
0&0&0&\cdots&1 &-3& 1&\cdots&0&0\\\hline
0&0&0&\cdots&0 &1& -2&\cdots&0&0\\
\vdots & \vdots & \vdots& \ddots&\vdots   & \vdots&  \vdots& \ddots& \vdots & \vdots  \\
0&0&0&\cdots&0 &0& 0&\cdots&-2&1\\
0&0&0&\cdots&0 &0& 0&\cdots&1&-2\\
\end{bmatrix}_{n\times n},$$ where $-3$ is at $(m,m)$ position $(m\leq n)$. Then $\det(R_{n,m})=(-1)^n(1+n+m+mn-m^2)$.
\end{lemma}
\noindent\textit{Proof.} Observe that $R_{n,m}=R_n-e_m{e_m}^T$. By Lemma \ref{Matdet}, we have \begin{eqnarray*}\det\left(R_{n,m}\right)&=&\det\left(R_n\right)-e_m^T\operatorname{adj}(R_n)e^m\\
&=&\det(R_n)-\mbox{ cofactor of the entry at } (m,m) \mbox{ position of } R_{n,m}\\
&=&\det(R_n)-(-1)^{m+m}\det\left(R_{m-1}\right)\cdot\det\left(R_{n-m}\right)\\
&=&(-1)^n(1+n)-(-1)^{m-1}m\cdot(-1)^{n-m}(1+n-m) \mbox{ (by Lemma \ref{det})}\\
&=&(-1)^n(1+n+m+mn-m^2).\end{eqnarray*}
\qed 

\begin{lemma}\label{vieta1} For $n\geq 2,$
$\sum\limits_{i=2}^{3n} \frac{1}{\rho_i}=-\frac{\alpha_{3n-2}}{\alpha_{3n-1}}=\frac{25n^2+30n-13}{60}\cdot$
\end{lemma}
\noindent\textit{Proof.} Let $M$ be a square matrix. By $M\{i\}$, we denote the submatrix of $M$, obtained by deleting the $i$th row and $i$th column of $M$. With this notation, we have from (\ref{char1}) that $\alpha_{3n-1}=\sum\limits_{i=1}^{3n} \det(-L_A\{i\}).$ 

Now, for $1\leq i \leq n$, $\det(-L_A\{i\})=\begin{vmatrix}
-2I_{n-1}& S \\ 
S^T & Q 
\end{vmatrix}$, where $S=-\sqrt{2}L_{V_0V_1}\{i\}$ and $Q= -L_{V_1V_1} - L_{V_1V_2}$. Using Schur complement, we have for $1\leq i \leq n$, $\det(-L_A\{i\})=\det\left(-2I_{n-1}\right)\cdot\det\left(Q+\frac{1}{2}S^TS\right)$.

Now,
$$Q+\frac{1}{2}S^TS= 
\begin{bmatrix}[ccccc|c|cccc]-2 &1&0&\cdots&0&0& 0&\cdots&0&1\\
1&-2 &1&\cdots&0&0& 0&\cdots&0&0\\
0 &1&-2&\cdots&0&0& 0&\cdots&0&0\\
\vdots & \vdots & \vdots& \ddots&\vdots   & \vdots&  \vdots& \ddots& \vdots & \vdots  \\
0&0&0&\cdots&-2 &1& 0&\cdots&0&0\\\hline 
0&0&0&\cdots&1 &-3& 1&\cdots&0&0\\\hline
0&0&0&\cdots&0 &1& -2&\cdots&0&0\\
\vdots & \vdots & \vdots& \ddots&\vdots   & \vdots&  \vdots& \ddots& \vdots & \vdots  \\
0&0&0&\cdots&0 &0& 0&\cdots&-2&1\\
1&0&0&\cdots&0 &0& 0&\cdots&1&-2\\
\end{bmatrix}_{2n\times 2n},$$ with $-3$ at $(2i,2i)$ position. 

Thus, $Q+\frac{1}{2}S^TS= R_{2n,2i}+e_1e_{2n}^T+e_{2n}e_1^T$. 

Let $R^1=R_{2n,2i}+e_1e_{2n}^T$. Then, 
\begin{eqnarray*}
\det\left(R^1\right)&=&\det\left(R_{2n,2i}\right)+e_{2n}^T\operatorname{adj}(R_{2n,2i})e_1\\
&=&\det\left(R_{2n,2i}\right)+\mbox{ cofactor of the entry at } (1,2n) \mbox{ position of } R_{2n,2i}\\
&=&\det(R_{2n,2i})+(-1)^{2n+1}\cdot 1\\
&=&1+2n+2i+4ni-4i^2 - 1 \mbox{ (by Lemma \ref{det1.1})}\\
&=&2n+2i+4ni-4i^2.\end{eqnarray*}

Thus, 
\begin{eqnarray*}
\det\left(Q+\frac{1}{2}S^TS\right)&=&\det\left(R^1\right)+e_{1}^T\operatorname{adj}(R^1)e_{2n}\\
&=&\det\left(R^1\right)+\mbox{ cofactor of the entry at } (2n,1) \mbox{ position of } R^1\\
&=&\det(R^1)+(-1)^{2n+1}\cdot \det
\begin{bmatrix}[cccccccccc]1&0&\cdots&0&0& 0&\cdots&0&1\\
-2 &1&\cdots&0&0& 0&\cdots&0&0\\
1&-2&\cdots&0&0& 0&\cdots&0&0\\
 \vdots & \vdots& \ddots&\vdots   & \vdots&  \vdots& \ddots& \vdots & \vdots  \\
0&0&\cdots&-2 &1& 0&\cdots&0&0\\
0&0&\cdots&1 &-3& 1&\cdots&0&0\\
0&0&\cdots&0 &1& -2&\cdots&0&0\\
 \vdots & \vdots& \ddots&\vdots   & \vdots&  \vdots& \ddots& \vdots & \vdots  \\
0&0&\cdots&0 &0& 0&\cdots&-2&1\\
\end{bmatrix}_{2n-1\times 2n-1}\\
&=&\det(R^1)-\left[1+\det(R_{2n-2, 2i-1})\right]\\
&=&2n+2i+4ni-4i^2 -\left[1+1+(2n-2)+(2i-1)+(4ni-2n-4i+2)\right.\\
&&\left.-(4i^2-4i+1)\right] \mbox{ (by Lemma \ref{det1.1})}\\
&=&2n.\end{eqnarray*}

Hence $\det(-L_A\{i\})=(-2)^{n-1}\cdot 2n=(-1)^{n-1}n\cdot 2^{n}$ for $i=1,2,\ldots,n.$

Again, for $n+1\leq i \leq 3n$, $\det(-L_A\{i\})=\begin{vmatrix}
-2I_{n}& U \\ 
U^T & Q\{i-n\} 
\end{vmatrix}$, where $U=-\sqrt{2}L_{V_0V_1}$ and $Q= -L_{V_1V_1} - L_{V_1V_2}$. Using Schur complement, we have for $n+1\leq i \leq 3n$,$\det(-L_A\{i\})=\det\left(-2I_n\right)\cdot\det\left(Q\{i-n\}+\frac{1}{2}U^TU\right)$.

Now, 
\begin{eqnarray*}Q\{i-n\}+\frac{1}{2}U^TU&=&
\begin{bmatrix}[ccccc|cccccc]-2 &1&0&\cdots&0&0& 0&\cdots&0&1\\
1&-2 &1&\cdots&0&0& 0&\cdots&0&0\\
0 &1&-2&\cdots&0&0& 0&\cdots&0&0\\
\vdots & \vdots & \vdots& \ddots&\vdots   & \vdots&  \vdots& \ddots& \vdots & \vdots  \\
0&0&0&\cdots&-2 &0& 0&\cdots&0&0\\\hline
0&0&0&\cdots&0 &-2& 1&\cdots&0&0\\
0&0&0&\cdots&0 &1& -2&\cdots&0&0\\
\vdots & \vdots & \vdots& \ddots&\vdots   & \vdots&  \vdots& \ddots& \vdots & \vdots  \\
0&0&0&\cdots&0 &0& 0&\cdots&-2&1\\
1&0&0&\cdots&0 &0& 0&\cdots&1&-2\\
\end{bmatrix}_{(2n-1)\times (2n-1)},\\
&&\mbox{\Big(the diagonal blocks are of size } (i-n-1)\times(i-n-1)\\
&& \mbox{ and } (3n-i)\times(3n-i) \mbox{ respectively.\Big)}\\
&=&R^2+e_1e_{2n-1}^T, \mbox{ where, } R^2=\begin{bmatrix}
R_{i-n-1} & \mathbf{0}\\
\mathbf{0} & R_{3n-i}
\end{bmatrix}+e_{2n-1}e_{1}^T.
\end{eqnarray*}

Clearly, $\det\left(R^2\right)
=\det(R_{n-i-1})\cdot\det(R_{3n-i})$ and hence,

\begin{eqnarray*}
\det\left(Q\{i-n\}+\frac{1}{2}U^TU\right)
&=&\det(R_{i-n-1})\det(R_{3n-i})+\mbox{ cofactor of the entry at the }\\
&&(1,2n-1) \mbox{ position of } R^2 \mbox{ (by Lemma \ref{Matdet})}\\
&=&\det(R_{i-n-1})\det(R_{3n-i})-\det(R_{i-n-2})\det(R_{3n-i-1})\\
&=&-(n-i)(3n-i+1)+(i-n-1)(3n-i)\\
&=&-2n.
\end{eqnarray*}

So, $\det(L_A\{i\})=(-1)^{n+1}n\cdot 2^{n+1}$  for $n+1\leq i \leq 3n$. 

Hence, 
\begin{eqnarray}
\nonumber     \alpha_{3n-1}&=&\sum\limits_{i=1}^{3n} \det(-L_A\{i\})\\
\nonumber     &=&\sum\limits_{i=1}^{n} \det(-L_A\{i\})+\sum\limits_{i=n+1}^{3n} \det(-L_A\{i\})\\
\nonumber     &=&n\cdot(-1)^{n+1}n\cdot 2^n+2n\cdot (-1)^{n+1}n\cdot 2^{n+1}\\
&=&(-1)^{n+1}2^n\cdot 5n^2.
\label{a3n-1}
\end{eqnarray}

Again suppose $M\{i,j\}$ denotes the principal submatrix of $M$ obtained by deleting the $i$th row and $j$th row and the corresponding columns. Then from (\ref{char1}), we have

\begin{eqnarray*}\alpha_{3n-2}&=&\sum\limits_{1\leq i <j\leq 3n} \det(-L_A\{i,j\})\\
&=&\sum\limits_{1\leq i <j\leq n} \det(-L_A\{i,j\})+\sum\limits_{n+1\leq i <j\leq 3n} \det(-L_A\{i,j\})+\sum\limits_{\substack{1\leq i <n\\ n+1\leq j\leq 3n}} \det(-L_A\{i,j\}).\end{eqnarray*}

For $1\leq i <j\leq n$, $$\det(-L_A\{i,j\})=\begin{vmatrix}
-2I_{n-2}& 0_{(n-2)\times 2n} \\ 
0_{2n\times (n-2)} & E_{2n\times 2n} 
\end{vmatrix},$$ where $$E=
\begin{bmatrix}[ccccc|c|ccc|c|cccc]-2 &1&0&\cdots&0&0& 0&\cdots&0&0& 0&\cdots&0&1\\
1&-2 &1&\cdots&0&0& 0&\cdots&0&0& 0&\cdots&0&0\\
0 &1&-2&\cdots&0&0& 0&\cdots&0&0& 0&\cdots&0&0\\
\vdots & \vdots & \vdots& \ddots&\vdots   & \vdots&  \vdots& \ddots&\vdots   & \vdots&  \vdots& \ddots& \vdots & \vdots  \\
0&0&0&\cdots&-2 &1& 0&\cdots&0&0& 0&\cdots&0&0\\\hline 
0&0&0&\cdots&1 &-3& 1&\cdots&0&\cdots& 0&\cdots&0&0\\\hline
0&0&0&\cdots&0 &1& -2&\cdots&0&0& 0&\cdots&0&0\\
\vdots & \vdots & \vdots& \ddots&\vdots   & \vdots&  \vdots& \ddots&\vdots   & \vdots&  \vdots& \ddots& \vdots & \vdots  \\
0&0&0&\cdots&0&0& 0&\cdots&-2 &1& 0&\cdots&0&0\\\hline 
0&0&0&\cdots&0&\cdots& 0&\cdots&1 &-3& 1&\cdots&0&0\\\hline
0&0&0&\cdots&0&0& 0&\cdots&0 &1& -2&\cdots&0&0\\
\vdots & \vdots & \vdots& \ddots&\vdots   & \vdots&  \vdots& \ddots&\vdots   & \vdots&  \vdots& \ddots& \vdots & \vdots  \\
0&0&0&\cdots&0 &0& 0&\cdots&0&0& 0&\cdots&-2&1\\
1&0&0&\cdots&0 &0& 0&\cdots&0&0& 0&\cdots&1&-2\\
\end{bmatrix}_{2n\times 2n},$$ with $-3$ at the $(2i, 2i)$ and at $(2j, 2j)$ positions.

By repeated application of Lemma \ref{Matdet}, we get 
$\det(E)=4n+8ij-4n(i-j)-4(i^2+j^2)$ and hence $\det(-L_A\{i,j\})=(-1)^{n-2}2^{n-2}[4n+8ij-4n(i-j)-4(i^2+j^2)]$ for $1\leq i <j\leq n$.

For $n+1\leq i <j\leq 3n$, let $p=i-n$ and $q=j-n$. So, $1\leq p <q\leq 2n$. 

Then, 
$\det(-L_A\{p,q\})=\begin{vmatrix}
-2I_{n}& \mathbf{0} \\ 
\mathbf{0} & F_{(2n-2)\times (2n-2)} 
\end{vmatrix},$  
where $F=\begin{pmatrix}
R_{p-1} & \mathbf{0} & e_{p-1}e_{2n-q-1}^T\\
\mathbf{0} & R_{q-p}&\mathbf{0}\\
e_{2n-q-1}e_{p-1}^T & \mathbf{0}& R_{2n-q-1}
\end{pmatrix}.$ 

Proceeding as before,  $\det(F)=2n(q-p)+2pq-(p^2+q^2)$ and hence $$\det(-L_A\{p,q\})=(-1)^n2^n[2n(q-p)+2pq-(p^2+q^2)] \mbox{ for } 1\leq p <q\leq 2n.$$
 
Similarly, for $1\leq i <n,$ $n+1\leq j\leq 3n,$ let $q=j-n$, i.e., for $1\leq q\leq 2n$,  

$$\det(-L_A\{i,q\})=\begin{vmatrix}
-2I_{n-1}& \mathbf{0} \\ 
\mathbf{0} & H_{(2n-1)\times (2n-1)} 
\end{vmatrix},$$ \\
where $\det(H)=\begin{cases}
			-q^2+4qi-4i^2+2n+2n(q-2i), & \text{if $2i\leq q$}\\
            	-q^2+4qi-4i^2+2n-2n(q-2i), & \text{if $2i> q$ }.
		 \end{cases}$

Hence,

\begin{eqnarray*}
\alpha_{3n-2}&=&(-1)^{n-2}2^{n-2}\frac{x^4+6x^3-7x^2}{3}+(-1)^{n}2^{n}\frac{4x^4-x^2}{3}+(-1)^n2^{n-1}\frac{4n^4+12n^3-n^2}{3}\\
&=&(-1)^n2^{n-2}\frac{25n^4+30n^3-13n^2}{3}\cdot
\end{eqnarray*}

Hence $\sum\limits_{i=2}^{3n} \frac{1}{\rho_i}=-\frac{\alpha_{3n-2}}{\alpha_{3n-1}}=\frac{(-1)^{n+1}2^{n-2}(25n^4+30n^3-13n^2)}{(-1)^{n+1}2^n\cdot 15n^2}=\frac{25n^2+30n-13}{60}\cdot$ \qed

\begin{lemma}\label{recurrence}
Let, $r_i$ be the determinant of $i\times i$ submatrix of $L^0_S$ formed by the first $i$ rows and first $i$ columns of $L^0_S$. Then for $1\leq i \le 2n$, $r_i=s_1(\sqrt{2}+\sqrt{3})^i+s_2(-\sqrt{2}-\sqrt{3})^i+s_3(\sqrt{2}-\sqrt{3})^i+s_4(-\sqrt{2}+\sqrt{3})^i$ where $s_1=\frac{(\sqrt{2}+\sqrt{3})(2+\sqrt{3})}{4\sqrt{6}}$, $s_2=\frac{(\sqrt{2}+\sqrt{3})(-2+\sqrt{3})}{4\sqrt{6}}$, $s_3=\frac{(\sqrt{2}-\sqrt{3})(-2+\sqrt{3})}{4\sqrt{6}}$ and $s_4=\frac{(\sqrt{2}-\sqrt{3})(2+\sqrt{3})}{4\sqrt{6}}\cdot$
\end{lemma}
\noindent\textit{Proof.} By observation $r_1=4$, $r_2=11$, $r_3=40$. If we choose $r_0=1$, then for $2\le i\le 2n$
\[r_{i}= 
     \begin{cases}
       3r_{i-1}-r_{i-2}, &\mbox{when } i\mbox{~is~even}\\
4r_{i-1}-r_{i-2}, &\mbox{when } i\mbox{~is~odd}.
     \end{cases}
\]

For $1\leq i \leq n-1$, let $c_i=r_{2i}$ and $d_i=r_{2i+1}$, then
\begin{eqnarray} \begin{cases}
       c_i=3d_{i-1}-c_{i-2}, &\mbox{when } i\mbox{~is~even} \label{case1}\\
d_i=4c_{i}-d_{i-1}, &\mbox{when } i\mbox{~is~odd}.
     \end{cases}
\end{eqnarray}
 
By substitution, we have

$$r_i=10r_{i-2}-r_{i-4},4\le i\le 2n.$$

The auxiliary equation for this recurrence relation is $x^4-10x^2+1=0.$ Solving, we get \begin{equation}
x=\pm (\sqrt{2}+\sqrt{3}), ~\pm (\sqrt{2}-\sqrt{3}).    
\end{equation}

Thus the general solution is $r_i=s_1(\sqrt{2}+\sqrt{3})^i+s_2(-\sqrt{2}-\sqrt{3})^i+s_3(\sqrt{2}-\sqrt{3})^i+s_4(-\sqrt{2}+\sqrt{3})^i.$
 Using the initial conditions $r_0=1$, $r_1=4$, $r_2=11$ and $r_3=40$ we get four linear equations in $s_1,s_2,s_3$ and $s_4$. Solving those equations by Cramer's rule, we get

\begin{eqnarray*}s_1&=&\frac{(\sqrt{2}+\sqrt{3})(2+\sqrt{3})}{4\sqrt{6}},\\ 
s_2&=&\frac{(\sqrt{2}+\sqrt{3})(-2+\sqrt{3})}{4\sqrt{6}},\\
s_3&=&\frac{(\sqrt{2}-\sqrt{3})(-2+\sqrt{3})}{4\sqrt{6}},\\
s_4&=&\frac{(\sqrt{2}-\sqrt{3})(2+\sqrt{3})}{4\sqrt{6}}\cdot
\end{eqnarray*}\qed

\begin{lemma}\label{vieta2} For $n\geq 2,$
 $\sum\limits_{i=2}^{3n} \frac{1}{\mu_i}=\frac{\beta_{2n-1}}{\det(L_S)}=\\ \frac{\frac{7\sqrt{6}}{192}\left[(4-2\sqrt{6})(5-2\sqrt{6})^{n-1}-(4+2\sqrt{6})(5+2\sqrt{6})^{n-1}\right]+\frac{7\sqrt{3}}{48}\left[(4n-2\sqrt{6}n+1)(\sqrt{3}-\sqrt{2})^{2n-1}+(4n+2\sqrt{6}n+1)(\sqrt{3}+\sqrt{2})^{2n-1}\right]}{(\sqrt{2}+\sqrt{3})^{2n}+(\sqrt{2}-\sqrt{3})^{2n}-2}$ 
and 

$\sum\limits_{i=2}^{3n} \frac{1}{\mu'_i}=\frac{\beta_{2n-1}'}{\det(L'_S)}=\\ \frac{\frac{7\sqrt{6}}{192}\left[(4-2\sqrt{6})(5-2\sqrt{6})^{n-1}-(4+2\sqrt{6})(5+2\sqrt{6})^{n-1}\right]+\frac{7\sqrt{3}}{48}\left[(4n-2\sqrt{6}n+1)(\sqrt{3}-\sqrt{2})^{2n-1}+(4n+2\sqrt{6}n+1)(\sqrt{3}+\sqrt{2})^{2n-1}\right]}{(\sqrt{2}+\sqrt{3})^{2n}+(\sqrt{2}-\sqrt{3})^{2n}+2}.$ 
\end{lemma}
\noindent\textit{Proof.} Apply Lemma \ref{Matdet} on $L_S$ and $L'_S$ to get, \begin{eqnarray}
     \det(L_s)&=&r_{2n}-r_{2n-2}-2\notag \\&=& (\sqrt{2}+\sqrt{3})^{2n}+(\sqrt{2}-\sqrt{3})^{2n}-2  
\label{det1}\end{eqnarray} 
and 
\begin{eqnarray}\det(L'_s)&=& r_{2n}-r_{2n-2}+2\notag \\&=&(\sqrt{2}+\sqrt{3})^{2n}+(\sqrt{2}-\sqrt{3})^{2n}+2\label{det2}\end{eqnarray} respectively.

It is not difficult to follow that
$\beta_{2n-1}=\beta'_{2n-1}=\sum\limits_{i=0}^{2n-1}r_i r'_{2n-1-i}$ where $r'_i$ is the determinant of $i\times i$ submatrix formed by last $i$ rows and last $i$ columns of $L^0_S$ for $1\leq i <2n.$ Note that ${r'_1}=4$, ${r'_2}=11$ and ${r'_3}=30.$ For convenience, we also choose $r'_0=1$.

Let $f(x)=\sum\limits_{i=0}^\infty r_ix_i$ and $g(x)=\sum\limits_{i=0}^\infty r'_ix_i.$ Then $\beta_{2n-1}=\mbox{coefficients of } x^{2n-1}-\mbox{coefficients of } x^{2n-3}$ in $f(x)g(x).$

So, \begin{eqnarray*}
       f(x)&=&\sum\limits_{i=0}^\infty r_ix_i=1+4x+11x^2+40x^3+\sum\limits_{i=4}^\infty r_ix_i\\
       &=& 1+4x+11x^2+40x^3+\sum\limits_{i=4}^\infty (10r_{i-2}-r_{i-4})x_i\\
       &=& 1+4x+11x^2+40x^3+10x^2\sum\limits_{i=4}^\infty r_{i-2}x^{i-2}-x^4\sum\limits_{i=4}^\infty r_{i-4}x^{i-4}\\
       &=& 1+4x+11x^2+40x^3+10x^2\big(f(x)-1-4x\big)-x^4f(x).\end{eqnarray*}
 
Hence, $f(x)=\frac{x^2+4x+1}{x^4-10x^2+1}\cdot$

In a similar way, we can deduce that $g(x)=\frac{x^2+3x+1}{x^4-10x^2+1}\cdot$

Thus we have \begin{eqnarray*}f(x)g(x)&=&\frac{x^4+7x^3+14x^2+14x+1}{(x^4-10x^2+1)^2}\cdot\end{eqnarray*}

Expressing as partial fraction,
\begin{eqnarray*}f(x)g(x)&=& \frac{\frac{12\sqrt{2}-7\sqrt{6}}{384}}{x+\sqrt{2}+\sqrt{3}}+\frac{\frac{-12\sqrt{2}-7\sqrt{6}}{384}}{x-\sqrt{2}-\sqrt{3}}+\frac{\frac{12\sqrt{2}+7\sqrt{6}}{384}}{x+\sqrt{2}-\sqrt{3}}+\frac{\frac{-12\sqrt{2}+7\sqrt{6}}{384}}{x+\sqrt{3}-\sqrt{2}}\\
&&+\frac{\frac{12-7\sqrt{3}}{192}}{(x+\sqrt{2}+\sqrt{3})^2}+\frac{\frac{12+7\sqrt{3}}{192}}{(x-\sqrt{2}-\sqrt{3})^2}+\frac{\frac{12+7\sqrt{3}}{192}}{(x-\sqrt{3}+\sqrt{2})^2}+\frac{\frac{12-7\sqrt{3}}{192}}{(x+\sqrt{3}-\sqrt{2})^2}.\end{eqnarray*}

Observe that the coefficient of $x^{2n-1}$ in\\
$\frac{1}{(x+\sqrt{2}+\sqrt{3})^2}=\frac{1}{(\sqrt{2}+\sqrt{3})^2}\left(1+\frac{x}{\sqrt{2}+\sqrt{3}}\right)^{-2}=(\sqrt{3}-\sqrt{2})^2\left(\sum\limits_{i=0}^\infty(\frac{-1}{\sqrt{2}+\sqrt{3}})^ix^i\right)^2$ is \\$(\sqrt{3}-\sqrt{2})^2\left(\sum\limits_{i=0}^{2n-1}(\frac{-1}{\sqrt{2}+\sqrt{3}})^i(\frac{-1}{\sqrt{2}+\sqrt{3}})^{2n-1-i}\right)=-2n(\sqrt{3}-\sqrt{2})^{2n+1}.$

Considering the coefficients of $x^{2n-1}$ and $x^{2n-3}$ from each fraction obtained in similar manner, we have $\beta_{2n-1}=\beta'_{2n-1}=\frac{7\sqrt{6}}{192}\left[(4-2\sqrt{6})(5-2\sqrt{6})^{n-1}-(4+2\sqrt{6})(5+2\sqrt{6})^{n-1}\right]$\\$+\frac{7\sqrt{3}}{48}\left[(4n-2\sqrt{6}n+1)(\sqrt{3}-\sqrt{2})^{2n-1}+(4n+2\sqrt{6}n+1)(\sqrt{3}+\sqrt{2})^{2n-1}\right].$

Hence $\sum\limits_{i=2}^{3n} \frac{1}{\mu_i}\\ =\frac{\frac{7\sqrt{6}}{192}\left[(4-2\sqrt{6})(5-2\sqrt{6})^{n-1}-(4+2\sqrt{6})(5+2\sqrt{6})^{n-1}\right]+\frac{7\sqrt{3}}{48}\left[(4n-2\sqrt{6}n+1)(\sqrt{3}-\sqrt{2})^{2n-1}+(4n+2\sqrt{6}n+1)(\sqrt{3}+\sqrt{2})^{2n-1}\right]}{(\sqrt{2}+\sqrt{3})^{2n}+(\sqrt{2}-\sqrt{3})^{2n}-2}$ 
and 

$\sum\limits_{i=2}^{3n} \frac{1}{\mu'_i}\\ =\frac{\frac{7\sqrt{6}}{192}\left[(4-2\sqrt{6})(5-2\sqrt{6})^{n-1}-(4+2\sqrt{6})(5+2\sqrt{6})^{n-1}\right]+\frac{7\sqrt{3}}{48}\left[(4n-2\sqrt{6}n+1)(\sqrt{3}-\sqrt{2})^{2n-1}+(4n+2\sqrt{6}n+1)(\sqrt{3}+\sqrt{2})^{2n-1}\right]}{(\sqrt{2}+\sqrt{3})^{2n}+(\sqrt{2}-\sqrt{3})^{2n}+2}.$ \qed

Combining Lemma \ref{vieta1} and Lemma \ref{vieta2}, we have the following theorem. 

\begin{theorem}\label{Kircchoff}
For $n \geq 2,$
\begin{eqnarray*}
Kf(P_n)&=&5n\left(-\frac{\alpha_{3n-2}}{\alpha_{3n-1}}+\frac{\beta_{2n-1}}{\det(L_S)}\right)\\
&=&5n\left(\frac{25n^2+30n-13}{60}+\frac{{{{\frac{7\sqrt{6}}{192}\left[(4-2\sqrt{6})(5-2\sqrt{6})^{n-1}-(4+2\sqrt{6})(5+2\sqrt{6})^{n-1}\right]}\atop{+\frac{7\sqrt{3}}{48}\left[(4n-2\sqrt{6}n+1)(\sqrt{3}-\sqrt{2})^{2n-1}+(4n+2\sqrt{6}n+1)(\sqrt{3}+\sqrt{2})^{2n-1}\right]}}}}{(\sqrt{2}+\sqrt{3})^{2n}+(\sqrt{2}-\sqrt{3})^{2n}-2}\right),    
\end{eqnarray*}

and 

\begin{eqnarray*}
Kf(P_n')&=&5n\left(-\frac{\alpha_{3n-2}}{\alpha_{3n-1}}+\frac{\beta'_{2n-1}}{\det(L_S')}\right)\\
&=&5n\left(\frac{25n^2+30n-13}{60}+\frac{{{{\frac{7\sqrt{6}}{192}\left[(4-2\sqrt{6})(5-2\sqrt{6})^{n-1}-(4+2\sqrt{6})(5+2\sqrt{6})^{n-1}\right]}\atop{+\frac{7\sqrt{3}}{48}\left[(4n-2\sqrt{6}n+1)(\sqrt{3}-\sqrt{2})^{2n-1}+(4n+2\sqrt{6}n+1)(\sqrt{3}+\sqrt{2})^{2n-1}\right]}}}}{(\sqrt{2}+\sqrt{3})^{2n}+(\sqrt{2}-\sqrt{3})^{2n}+2}\right)\cdot    
\end{eqnarray*}

\end{theorem}

In $1997$, Chung \cite{Chu97} showed that the number of spanning trees $\tau(G)$ for a connected graph $G$ is equal to the product of non-zero Laplacian eigenvalues of $G$ devided by the number of vertices of $G$. From (\ref{a3n-1}), (\ref{det1}) and (\ref{det2}) we have the following theorem.
\begin{theorem}
The total number of spanning trees of $P_n$ is 
\begin{eqnarray*}\tau(P_n)&=&\frac{\prod\limits_{i=2}^{3n} \lambda_i \prod\limits_{j=1}^{2n}\mu_j}{5n}\\
&=&\frac{(-1)^{3n-1}\alpha_{3n-1}\det{(L_s)}}{5n}\\
&=&2^n n\left((\sqrt{2}+\sqrt{3})^{2n}+(\sqrt{2}-\sqrt{3})^{2n}-2\right)
\end{eqnarray*}

and the total number of spanning trees of $P'_n$ is \begin{eqnarray*}\tau(P'_n)&=&\frac{\prod\limits_{i=2}^{3n} \lambda_i \prod\limits_{j=1}^{2n}\mu'_j}{5n}\\
&=&\frac{(-1)^{3n-1}\alpha_{3n-1}\det{(L'_s)}}{5n}\\
&=&2^n n\left((\sqrt{2}+\sqrt{3})^{2n}+(\sqrt{2}-\sqrt{3})^{2n}+2\right).
\end{eqnarray*}
\end{theorem}

\section{Relation between Kirchhoff index and Wiener index}
First we calculate Wiener index for $P_n$ and $P'_n.$

\begin{theorem}\label{Wiener}
The Wiener index of $P_n$ is
\[ W(P_n)=\begin{cases} 
      \frac{25}{4}n^3+9n^2, \textit{~~when~} n\textit{~is~even}\\
      \frac{25}{4}n^3+9n^2-\frac{n}{4}, \textit{~~when~} n\textit{~is~odd.} 
   \end{cases}
\]
and the Wiener index of $P'_n$ is
\[ W(P'_n)=\begin{cases} 
      \frac{25}{4}n^3+9n^2-2n, \textit{~~when~} n\textit{~is~even}\\
      \frac{25}{4}n^3+9n^2-\frac{9n}{4}, \textit{~~when~} n\textit{~is~odd.} 
   \end{cases}
\]
\end{theorem}
\noindent\textit{Proof.} 
There are three type of vertices of $P_n,$
\begin{enumerate}
  \item $a$-type: vertex with degree $3$ and non-adjacent to vertex of degree $2$,
  \item $b$-type: vertex with degree $3$ and which is adjacent to a vertex of degree $2$,
  \item $c$-type: vertex of degree $2$ (middle vertex) of $P_n$.
\end{enumerate}

Now we observe the following.

\begin{enumerate}
   \item Sum of the distances from an $a$-type vertex to
  
   \begin{enumerate}
     \item all the upper vertices $j:2\sum\limits_{i=1}^{n-1} i+n=n^2.$
     \item all the lower vertices $j':2\sum\limits_{i=1}^{n} i+n=n^2+2n.$
     \item all the middle vertices $\tilde{j}$:
     \begin{enumerate}
       \item when $n$ is even: $4\sum\limits_{i=1}^{n/2} i=\frac{n}{2}(n+2).$ 
       \item when $n$ is odd: $4\sum\limits_{i=1}^{(n+1)/2} i-(n+1)=\frac{1}{2}(n+1)^2.$
     \end{enumerate}
   \end{enumerate}
 
   \item Sum of the distances from a $b$-type vertex to
  
   \begin{enumerate}
     \item all the upper vertices $j:2\sum\limits_{i=1}^{n-1} i+n=n^2.$
     \item all the lower vertices  $j':2\sum\limits_{i=1}^{n} i+(n+1)=(n+1)^2.$
     \item all the middle vertices $\tilde{j}$:
     \begin{enumerate}
       \item when $n$ is even: $2\sum\limits_{i=1}^{n/2} (2i-1)+n=\frac{n^2}{2}+n.$ 
       \item when $n$ is odd: $2\sum\limits_{i=1}^{(n+1)/2} (2i-1)-1=\frac{(n+1)^2-2}{2}+n.$ 
     \end{enumerate}
   \end{enumerate}

 \item Sum of the distances from a $c$-type vertex to
  
   \begin{enumerate}
     \item all the upper vertices $j:2\sum\limits_{i=1}^{n} i+n=n^2+2n.$
     \item all the lower vertices $j':2\sum\limits_{i=1}^{n} i+n=n^2+2n.$
     \item all the middle vertices $\tilde{j}$:
     \begin{enumerate}
       \item when $n$ is even: $4\sum\limits_{i=1}^{n/2} i+n-2=\frac{n^2+4n-4}{2}.$ 
       \item when $n$ is odd: $4\sum\limits_{i=1}^{(n+1)/2} i-4=\frac{(n+1)(n+3)}{2}-4.$  
     \end{enumerate}
   \end{enumerate}

 \end{enumerate}
Hence the Wiener index for $P_n$ is 
\begin{eqnarray*}
W(P_n)&=&\frac{\sum\limits_{j=1}^{2n}\left(2n^2+ (n^2+2n)+(n+1)^2+\frac{n}{2}(n+2)+\frac{n^2}{2}+n+n^2+2n   \right)+\sum\limits_{j=1}^{2n}\left(\frac{n^2+4n-4}{2}\right)    }{2}
\\&=&\frac{2n\left(2n^2+ (n^2+2n)+(n+1)^2+\frac{n}{2}(n+2)+\frac{n^2}{2}+n+n^2+2n   \right)+n\left(\frac{n^2+4n-4}{2}\right)    }{2}\\
&=&\frac{25}{4}n^3+9n^2,
\end{eqnarray*} when $n$ is even
and

\begin{eqnarray*}
W(P_n)&=&\frac{\sum\limits_{j=1}^{2n}\left(2n^2+ (n^2+2n)+\frac{1}{2}(n+1)^2+(n+1)^2+\frac{(n+1)^2-2}{2}+n^2+2n   \right)+\sum\limits_{j=1}^{n}\left(\frac{(n+1)(n+3)}{2}-4\right)    }{2}\\
&=&\frac{2n\left(2n^2+ (n^2+2n)+\frac{1}{2}(n+1)^2+(n+1)^2+\frac{(n+1)^2-2}{2}+n^2+2n   \right)+n\left(\frac{(n+1)(n+3)}{2}-4\right)    }{2}\\
&=&\frac{25}{4}n^3+9n^2-\frac{n}{4},\end{eqnarray*} when $n$ is odd. 

By similar approach we get the Wiener index for $P'_n$ as $$W(P'_n)=\left\{\begin{array}{ll}\frac{25}{4}n^3+9n^2-2n,& \mbox{when } n \mbox{ is even,}\\
\frac{25}{4}n^3+9n^2-\frac{9n}{4},& \mbox{when } n \mbox{ is odd}.
\end{array}
\right.$$ \qed

From Theorem $\ref{Kircchoff}$ and Theorem $\ref{Wiener}$ we have the following theorem. 
\begin{theorem}
$$\lim_{n \to\infty}\frac{W(P_n)}{Kf(P_n)}=3$$
and $$\lim_{n \to\infty}\frac{W(P'_n)}{Kf(P'_n)}=3.$$
\end{theorem}

In Table \ref{tab2}, for different values of $n$, we list the values of Kirchhoff indices and Wiener indices and the ratios of Wiener indices to Kirchhoff indices for both $P_n$ and $P_n'$. It is evident from the table that both the ratios are gradually approaching to 3.

\begin{table}[h]
\begin{center}
\begin{tabular}{c|c|c|c|c|c|c}
\hline
 $n$ &  $Kf(P_n)$ & $W(P_n)$ & $\frac{W(P_n)}{Kf(P_n)}$   &  $Kf(P'_n)$ & $W(P'_n)$ & $\frac{W(P'_n)}{Kf(P'_n)}$\\ 
 \hline
 2&             39.083333 &86 &2.200426458       & 38.5 & 82 &2.12987012987 \\
 \hline
 3  & 107.715909 &249 &2.3116362505    & 107.583333  & 243 & 2.25871418206\\
 \hline
4 & 226.166667 &544 &2.40530581812      &  226.142857 &536 & 2.37018319796\\
 \hline
5 & 406.806193 &1005 & 2.47046386533  &  406.802434 &995 &2.44590473616 \\
 \hline
6 & 662.098485 & 1674  &2.52832477029      &  662.097938 &1662 & 2.51020265222 \\
\hline
7 & 1004.536492 &2583 & 2.57133515862  &  1004.536417 & 2569 & 2.55739857363 \\
 \hline
8 & 1446.619048 &3776 & 2.61022416732        &  1446.619038 &3760 & 2.5991639134\\
 \hline
9 & 2000.845977 &5283 & 2.64038314829 &  2000.845975 &5265 &2.63138695621 \\
 \hline
10 & 2679.717254 &7150 & 2.66819194799       &  2679.717254 &7130 & 2.660728474 \\
 \hline
20 & 19073.869017 &53600 & 2.81012729783        &  19073.869017 & 53560 & 2.80803018791\\
\hline
99 & 2080862.36308 &6152553 & 2.95673231885        &  2080862.36308 & 6152355 & 2.95663768188\\
 \end{tabular}
\caption{\label{tab2}Kirchhoff index and Wiener index of pentagonal cylinder chain $P_n$ and pentagonal M\"{o}bius chain $P'_n$ for different values of $n\geq 2$.}
\end{center}
\end{table}

\section{Concluding remarks}
In this paper, we have derived explicit formulas for Kirchhoff index and Wiener index of linear pentagonal cylinder chain graph $P_n$ of $5n$ vertices and linear pentagonal M\"{o}bius chain graph $P'_n$ of $5n$ vertices. We have also established that for large values of $n$, Wiener index is almost three times the Kirchhoff index for both the graphs. 

\section*{Acknowledgement} University Grants Commission, India has provided a partial support through Senior Research Fellowship to the first author for this work.

\end{document}